\documentclass[11pt]{article}
\usepackage{amsmath}
\usepackage{graphicx}
\usepackage{latexsym}
\usepackage{amssymb}
\usepackage{color}

\newlength{\titleright}
\def\be{\begin{equation}}
\def\ee{\end{equation}}
\def\bea{\begin{eqnarray}}
\def\ba{\begin{array}{l}\displaystyle}
\def\eea{\end{eqnarray}}
\def\ea{\end{array}}

\def\frac#1#2{{#1 \over #2}}

\def\r{\rho}
\def\e{\varepsilon}
\def\l{\lambda}
\def\bx{{\bf x}}
\def\bv{{\bf v}}
\def\bw{{\bf w}}
\def\tH{{\widetilde H}}
\def\phi{{\varphi}}
\def\vpp{{\varPhi}}
\def\RR{\mathbb R}
\def\del{{\partial}}
\newcommand{\wh}{\widehat}
\newcommand{\wt}{\widetilde}
\newtheorem{definition}{Definition}[section]

\numberwithin{equation}{section}

\begin{document}
\title{Implicit-Explicit Runge-Kutta schemes for the Boltzmann-Poisson system for semiconductors}

\author{
Giacomo Dimarco\footnote{ 
Institut de Math\'{e}matiques de Toulouse, Toulouse, France. 
\texttt{giacomo.dimarco@math.univ-toulouse.fr}},
Lorenzo Pareschi\footnote{Department of Mathematics and Computer Science, University of Ferrara, via Ma\-chia\-vel\-li  35, 44121 Ferrara, Italy.
\texttt{lorenzo.pareschi@unife.it}},
Vittorio Rispoli\footnote{Department of Mathematics and Computer Science, University of Ferrara, via Ma\-chia\-vel\-li  35, 44121 Ferrara, Italy.
\texttt{rspvtr@unife.it}}
}
\maketitle

\begin{abstract}
In this paper we develop a class of Implicit-Explicit Runge-Kutta schemes for solving the multi-scale semiconductor Boltzmann equation.
The relevant scale which characterizes this kind of problems is the diffusive scaling.
This means that, in the limit of zero mean free path, the system is governed by a drift-diffusion equation.
Our aim is to develop a method which accurately works for the different regimes encountered in general semiconductor simulations:
the kinetic, the intermediate and the diffusive one.
Moreover, we want to overcome the restrictive time step conditions of standard time integration techniques
when applied to the solution of this kind of phenomena without any deterioration in the accuracy.
As a result, we obtain high order time and space discretization schemes which do not suffer from the usual parabolic stiffness in the diffusive limit.
We show different numerical results which permit to appreciate the performances of the proposed schemes.
\end{abstract}


\maketitle

\section{Introduction}\label{sec1}
The application of kinetic theory to the modeling of semiconductor devices simulations is a long dated but still very
active research field because of the richness and diversity of observed phenomena \cite{Mar, MRS}.
From the mathematical point of view, semiconductor devices can be accurately described by kinetic equations when the mean free path
of the particles is large compared to a macroscopic characteristic length of the system (\cite{ARR, CGMS2, CGMS1}).
On the other hand, when the average time between particles' collisions is small, the relevant scaling is the diffusive one and in such regimes
these systems can be described by macroscopic drift-diffusion equations \cite{Pou, MRS, SZ, FJO}.
Unfortunately, this passage from the microscopic to the macroscopic description leads to challenging numerical difficulties.
Indeed, when the macroscopic reference length is several orders of magnitude larger than the mean free path, the kinetic equation contains
stiff terms. This means that classical numerical methods need time step restrictions which make their use prohibitively expensive.
In these cases, it becomes attractive to use domain decomposition strategies, which are able to solve the microscopic and the macroscopic
models wherever it is necessary.
These approaches have been largely studied for kinetic equations both for the diffusive \cite{BAL,DEG1,DEG3,Klar} and for the hydrodynamic
scaling \cite{DEG2, DIM1, DIM2}.
However, even if these methods are very efficient they are affected by some difficulties due to the fact that it is not always a simple task
to define the different regions of the domain in which the use of a macroscopic model is fully justified.
Thus in the recent past alternative strategies have been studied, the so-called Asymptotic Preserving (AP) schemes.
They consist in solving the original kinetic model in the full domain avoiding the time step restriction caused by the presence of different stiff terms in the equations.
The AP methods automatically transform the original problem in the numerical approximations of the relevant macroscopic model when the
scaling parameter goes to zero \cite{J, FilbJin, JP, JPT1, JPT2, Kl2, BEN, LemMieu2008, BC, DP, GOS1, GOS2, NP1, NP2}.
Recently, this approach has been considered in the framework of Implicit-Explicit (IMEX) Runge-Kutta schemes with the aim of deriving
high order numerical methods which are accurate in all regimes \cite{ARS, CK, DP12, DP1, BPR12, PR05}.
This means that such schemes are able to preserve the desired order of accuracy even in the limit when the scaling parameter tends to zero.

In this paper we develop high order schemes for the resolution of the Boltzmann-Poisson semiconductor model.
In this kind of problems, since the characteristic speed of the hyperbolic part of the kinetic equation is of the order of $1/\e$
(where $\e$ is proportional to the mean free path), the CFL condition for an explicit approach would require $\Delta t = O(\e\,\Delta x)$.
Of course, in the diffusive regime where $\e\ll\Delta x$, this would be too much restrictive since a parabolic condition
$\Delta t = O(\Delta x^2)$ would suffice to solve the limiting equation. This kind of difficulties has been studied in \cite{JP, JPT1, JPT2, Kl2, NP1, NP2} using different semi-implicit approaches. However, most of the previous literature on the subject, originates consistent low order explicit schemes for the limit model. Such explicit schemes clearly suffer from the usual stability restriction $\Delta t = O(\Delta x^2)$. Here, on the contrary to previous approaches, we are able to guarantee for all regimes a linear time step limitation $\Delta t = O(\Delta x)$, independently from $\e$, and high order accuracy in time and space. In order to accomplish this task, we first rewrite our system using the parity formalism then, following \cite{BPR12}, we add and subtract the limiting diffusive flux to the convective kinetic flux. We then discretize the reformulated problem by Implicit-Explicit Runge Kutta method. As we will show this will allow to get Asymptotic Preserving high order schemes which uniformly work for the
different values of $\e$ and which automatically originate an IMEX Runge-Kutta method for the limiting convection-diffusion equation in which the diffusive term is discretized implicitly. Finally, in order to permit realistic simulations we consider the challenging case of complex interaction operators. These terms are handled by introducing a penalization technique which permit to avoid the inversion of such operators when their stiff character suggests an implicit treatment.

The rest of the paper is organized as follows.
First we present the kinetic semiconductor equation, its drift-diffusion limit and the reformulated
parity system.
The third section provides the basis of the numerical method describing in details the time discretization.
The fourth section describes phase-space variables discretization, based on conservative finite difference schemes for
the space variables and on a Gauss-Hermite approximation for the velocity variables.
Numerical results for the proposed schemes are presented in the fifth section. A concluding section ends the paper.

\section{The Boltzmann equation and its drift-diffusion limit}
We consider the Boltzmann equation under the diffusive scaling which describes the time evolution of electrons inside semiconductor devices.
Let $f(t, \bx, \bv)$ be the density distribution function for particles at time $t \geq 0$, where position and velocity variables $\bx$ and
$\bv$ are such that $(\bx,\bv)\in\Omega\times \RR^d$, with $\Omega \subset \RR^d$ and $d = 1, 2,$ or $3$. Under these assumptions, the time
evolution of the system is described by \cite{Klar}
\begin{equation}\label{tone}
\e\,\partial_t f +\bv\cdot\nabla_\bx f - \frac{q}{m}E \cdot \nabla_\bv f = {1\over\e}Q(f) + \e \widehat G(f) \,.
\end{equation}

In this formula, $\widehat G = \widehat G(f)$ is an integrable function of $\bv$ which models the generation and recombination process inside the semiconductor,
$\e$ is proportional to the mean free path and $E(t,\bx) = -\nabla_\bx \vpp(t,\bx)$ is the electric field which is self-consistently computed solving a suitable Poisson
equation for the electric potential $\vpp$.
Constants $q$ and $m$ are respectively the elementary charge and the effective mass of the electrons.
The anisotropic collision term $Q(f)$ is defined by
\be \label{eq:Q}
  Q(f)(\bv) = \int \sigma(\bv, \bw) \{M(\bv) f (\bw)- M(\bw) f(\bv)\} \,d\bw,
\ee
where $M$ is the constant in time normalized Maxwellian at temperature $\theta$
$$  M(\bv) = \frac1{(2\pi\theta)^{d/2}} \exp \left( - {{|\bv|^2}\over{2\theta}}\right)  $$
and $\sigma$ is the anisotropic scattering kernel which is rotationally invariant and satisfies
$$  \sigma( \bv, \bw ) = \sigma( \bw, \bv ) \geq s_0 > 0  $$
for some given constant $s_0$.
We also assume that the collision frequency $\l$ satisfies the following bound for some positive constant $\l_M$ \cite{Kl2}
\begin{equation}\label{collf}
0 < \, \l (\bv) = \int \sigma( \bv, \bw ) \, M(\bw)\, d\bw \, \leq \l_M.
\end{equation}
The collision operator is bounded and nonnegative on the Hilbert space, see \cite{Pou,Kl2},\\
$H=\mathcal{L}^2(\RR^d,M^{-1}(\bv)d\bv)$ and it has a one-dimensional kernel spanned by $M$.
Moreover, if we let $h_i \in H$ be the unique solutions to the problem
$$  Q(h_i\,M) = - v_i M \,, \qquad \int h_i(\bv)\,M(\bv)\, d\bv = 0, \quad i=1,\ldots,d\, ,$$
then, since $\sigma$ is rotationally invariant, it follows that there exists a positive constant $D$ such that
\begin{equation}\label{diffc}
  \int v_i h_j M(\bv) \, d\bv = D\delta_{ij}, \qquad i,j=1,\ldots,d\,.
\end{equation}
Thus, defining the total mass $\r=\rho(t,\bx)$ as
$$ \rho=\int f(\bv)\,d\bv\,, $$
one can show (using for example the Hilbert expansion $f = f_0 +\e f_1+\e^2 f_2$ with $Q(f_0)=0$, see \cite{LemMieu2008} for instance)
that when $\e \to 0$, $f(t, \bx, \bv)$ is approximated by $\rho(t, \bx) M(\bv)$ with $\r$ satisfying the drift-diffusion
equation \cite{Klar, MRS}
\begin{equation}\label{televen}
 \partial_t\rho= \nabla_\bx \cdot \left ( D \nabla_\bx \rho  + \eta \rho E \right) + \widetilde{G}.
\end{equation}
In this equation, $D$ is the diffusion coefficient defined implicitly in terms of the cross section by (\ref{diffc}), the constant $\eta$
is the so-called mobility given by the Einstein relation $q\,D=\eta\,m\,\theta$ and $\widetilde{G}$ is the integral of the generation
recombination function
$$  \wt{G}(t,\bx)= \int \widehat G(f(\bv))\, d \bv. $$

\subsection{Even and odd parities}
We now define the even and odd parities formalism which we will use in the development of the numerical schemes.
To this aim we split equation (\ref{tone}) into two equations, one for $\bv$ and one for~$-\bv$
\begin{equation}\label{thtwo}
\begin{split}
\e\, \partial_t f +\bv\cdot\nabla_\bx f - {q\over m}E \cdot \nabla_\bv f & = {1\over\e}Q(f)(\bv) + \e \widehat G(f),\\
\e\, \partial_t f -\bv\cdot\nabla_\bx f + {q\over m}E \cdot \nabla_\bv f & = {1\over\e}Q(f)(-\bv) + \e \widehat G(f).
\end{split}
\end{equation}
Next, we introduce the so called even parity $r$ and odd parity $j$ defined by
\begin{align}
r(t,\bx, \bv) & = \frac12 \, \Big( f(t,\bx, \bv)+f(t,\bx,-\bv) \Big),       \label{eq:parir}  \\
j(t,\bx, \bv) & = \frac1{2\e} \, \Big( f(t,\bx, \bv)-f(t, \bx,-\bv) \Big).  \label{eq:parij}
\end{align}
Now, adding and subtracting the two equations in (\ref{thtwo}) leads to
\begin{equation}\label{thsix}
\begin{split}
  \partial_t r + \bv \cdot \nabla_\bx j - {q\over m} E \cdot \nabla_\bv j & = {1\over\e^2} Q(r) + \widehat G(r),\\
  \partial_t j+ { 1\over\e^2} \left( \bv \cdot \nabla_\bx r - {q\over m}E \cdot \nabla_\bv r \right) & = - {1\over\e^2} \, \l j \,,
\end{split}
\end{equation}
where $\lambda$ is the collision frequency defined in (\ref{collf}) and where we used the property that
$$\int \sigma(\bv, \bw) j(\bw)\,d \bw=0$$
since $j$ is an odd function. We also assume $\widehat{G}=\wh{G}(r)$ is an even function, which is a physically consistent choice in semiconductor
simulations. An important advantage of this formulation is that now only one time scale appears in our new system (\ref{thsix}).

From the above formulation it is easy to see that in the limit $\e\rightarrow 0$, formally we get that the integral
of the even parity $r$ with respect to the velocity variable is the solution to (\ref{televen}).
In fact, the macroscopic variable $\rho$ can be expressed in terms of $r$ as
\be
\rho=\int f(\bv)d\bv=\int r(\bv)d\bv.
\ee
Taking the formal limit $\e = 0$ in system (\ref{thsix}) we get
\begin{eqnarray}
Q(r)  & = &  0 \,, \label{thq:1}\\
\l j    & = &  - \bv \cdot \nabla_\bx r + {q\over m} E \cdot \nabla_\bv r \,. \label{thq:2}
\end{eqnarray}
The first equation implies that $r=r(t,\bx,\bv)=\r(t,\bx) M(\bv)$ being Ker$(Q)=$ Span$\{M\}$.
Now, replacing equation (\ref{thq:2}) and $r=\r M$ in the first equation of system (\ref{thsix}) and integrating
over the velocity space gives the drift diffusion equation (\ref{televen}).
From such computation we recover, moreover, that $D$ is given by the following formula
$$
  D = \int \frac{|\bv|^2M}{\l}\,d\bv.
$$
In the following, in order to simplify notations, we will assume $q=m=1$ and $2\theta=1$, which implies $D=\theta\eta=\eta/2$.

\section{Time discretization of the multi-scale system}
We recall here the system (\ref{thsix}) that we want to solve
\begin{eqnarray}
  \partial_t r + \bv \cdot \nabla_\bx j - E \cdot \nabla_\bv j & = & {1\over\e^2} Q(r) + \widehat G(r)        \label{sys_t_solv} \\
  \partial_t j+ { 1\over\e^2} ( \bv \cdot \nabla_\bx r - E \cdot \nabla_\bv r) & = & - {1\over\e^2} \, \l \, j  \label{sys_t_solv2}
\end{eqnarray}
and summarize the properties we demand to our time discretization scheme:
\begin{enumerate}
 \item  The scheme has to be asymptotic preserving (AP). This ensures stability condition independently from $\e$. In other words, the scheme should satisfy the following definition
  \begin{definition}
   A consistent time discretization method for (\ref{sys_t_solv}-\ref{sys_t_solv2}) of stepsize $\Delta t$ is {\em asymptotic
   preserving (AP)} if, independently of the initial data and of the stepsize $\Delta t$, in the limit $\varepsilon\to 0$ becomes
   a consistent time discretization for equation (\ref{televen}).
  \end{definition}
 \item The scheme has to be an high order asymptotically accurate (AA) method. This means we want our scheme to satisfy the following definition
  \begin{definition}
    A consistent time discretization method for (\ref{sys_t_solv}-\ref{sys_t_solv2})) of stepsize $\Delta t$ is {\em asymptotically accurate (AA)} if, independently of
   the initial data and of the stepsize $\Delta t$, the order of accuracy in time is preserved for $\e\rightarrow 0$.
  \end{definition}
 \item The scheme should solve, in the limit $\e \rightarrow 0$, the drift-diffusion equation (\ref{televen}) with an implicit treatment of the diffusion term.
 This ensures a stability condition for the time step which is of the order of the space discretization: $\Delta t = O(\Delta x)$.
 \item We want to avoid the difficult inversion of complex collision operators which occurs when classical implicit solvers are used.
\end{enumerate}

\subsection{The Boscarino-Pareschi-Russo reformulation}

We now focus on the first three points of the above list while we will discuss the last point at the end of the section.
In order to construct time discretizations which satisfy the above requirements we reformulate system (\ref{sys_t_solv}-\ref{sys_t_solv2}) following the approach introduced in \cite{BPR12}. We add to both sides of equation
(\ref{sys_t_solv}) a term that permits to choose which kind of time integrator (if explicit or implicit) we intend to use for the diffusive term in the drift-diffusion equation (\ref{televen}). This term reads
$$
\bv\cdot\nabla_\bx\left(\mu{\bv\over\l}\cdot\nabla_\bx r\right),
$$
where $\mu = \mu(\e)$ is a positive function such that $\mu(0)=1$.
The modified system reads
 \begin{equation}\label{IEsys}
  \begin{split}
   \partial_t r & + \bv \cdot \nabla_\bx \left( j + \mu{\bv\over\l}\cdot\nabla_\bx r\right) - E \cdot \nabla_\bv j
           = {1\over\e^2} Q(r) + \bv\cdot\nabla_\bx \left(\mu{\bv\over\l}\cdot\nabla_\bx r\right) + \widehat G(r),  \\
   \partial_t j & = - { 1\over\e^2} \Big( \bv \cdot \nabla_\bx r - E \cdot \nabla_\bv r \Big) - {1\over\e^2} \, \l j \,.
  \end{split}
 \end{equation}
As we will see, the introduction of this new term allows to avoid the parabolic time step limitations for the limit drift diffusion
equation (\ref{televen}). We will discuss later the different possible choices for $\mu$.

In order to write the general time discretization formulation, we rewrite the previous system (\ref{IEsys}) in a more compact form as
\begin{equation}\label{syssimp}
 \begin{split}
  \partial_t r & = f_1(r,j) + {1\over\e^2} Q(r) + f_2(r),\\
  \partial_t j & = -{ 1\over\e^2} g(r,j)
 \end{split}
\end{equation}
where we defined
\begin{eqnarray*}
 f_1(r,j) & = & -\bv \cdot \nabla_\bx \left( j + \mu {\bv\over\l} \cdot \nabla_\bx r \right) + E \cdot \nabla_\bv j + \widehat G(r), \\
 f_2(r)   & = &  \mu \, \bv \cdot \nabla_\bx \left( {\bv\over\l} \cdot \nabla_\bx r \right), \\
 g(r,j)   & = &   \l j + \left( \bv \cdot \nabla_\bx r - E \cdot \nabla_\bv r\right).
\end{eqnarray*}

\noindent {\bf Remark:} We point out that in this paper we consider the case in which $\mu$ does not depend on the space variable.
More precisely, $\mu$ could depend on $\e$ but not on $x$ and $t$. Since in the general case $\e$ instead depends on $x$ and $t$,
in the work we are assuming that in the space-time domain the range in which $\e$ may vary is small enough to
allow us to keep $\mu$ constant. We numerically observed that this condition is, indeed, true for a large set of values of $\e$.
We remind to a future work the development of schemes which could take into account the possibility for $\mu$ to vary.

\subsection{IMEX Runge-Kutta schemes}

An IMEX Runge-Kutta scheme \cite{ARS, BPR12} applied to the above system reads for the internal stages $k=1,\ldots,\nu$ as
\begin{eqnarray}
 R^{(k)} & = & r^n + \Delta t \sum_{j=1}^{k-1} \widetilde{a}_{kj}\,f_1\left(R^{(j)}, J^{(j)}\right)
           +   \Delta t \sum_{j=1}^k a_{kj} \left( {1 \over\e^2} Q \left( R^{(j)} \right)
           +   f_2 \left( R^{(j)} \right) \right)    \label{stRK} \\
 J^{(k)} & = & j^n - {\Delta t\over\e^2} \sum_{j=1}^k a_{kj} \, g \left( R^{(j)}, J^{(j)} \right)  \label{stJK}
 \end{eqnarray}
while the numerical solution is given by
\begin{eqnarray}
 r^{n+1} & = & r^n + \Delta t\sum_{k=1}^\nu\widetilde w_k f_1\left(R^{(k)}, J^{(k)}\right) +
               \Delta t \sum_{k=1}^\nu w_k \left( {1\over\e^2}Q\left(R^{(k)}\right) + f_2\left(R^{(k)}\right) \right)  \label{numsolr}  \\
 j^{n+1} & = & j^n - {\Delta t\over\e^2} \sum_{k=1}^\nu w_k\, g\left(R^{(k)}, J^{(k)}\right).  \label{numsolj}
\end{eqnarray}
In the above formulas, matrices $ \wt A=(\wt a_{ij} )$, $\wt a_{ij} = 0$ for $j \geq i$ and $A = (a_{ij})$ are $\nu\times \nu$ matrices such
that the resulting scheme is explicit in $f_1$ and implicit in $Q(r)$, $f_2$ and $g$.
In general, an IMEX Runge-Kutta scheme, is characterized by the above defined two matrices and by the coefficient vectors
$\wt w =( \widetilde{w}_{1},..,\wt w_{\nu})^{T}$, $w =(w_{1},..,w_{\nu})^{T}$.
Since computational efficiency is of paramount importance, in the sequel we restrict our analysis to diagonally implicit Runge-Kutta
(DIRK) schemes for the source terms ($a_{ij} = 0,$ for $j>i$).
The use of a DIRK scheme is enough to ensure that the two transport terms in the two equations of system (\ref{IEsys})
are always explicitly evaluated.
In fact, observe that, when the integration for the odd parity $j$ is performed the values of $r$ are already available for the corresponding
stage.
This is due to the use of a partitioned Runge-Kutta approach for the time integration of our system.
The type of schemes introduced can be represented with a compact notation by a double Butcher tableau,
$$
  \begin{tabular}{l|r c}
   $\wt c$  &  $\wt A$  \\
   \hline
            &  $\wt w^T$
  \end{tabular} \ \ \ \ \ \ \qquad
  \begin{tabular}{l|r c}
   $c$      &  $A$      \\
   \hline
            &  $w^{T}$
  \end{tabular}
$$
where coefficients $\wt c$ and $c$ are given by the usual relation
$\wt c_i = \sum_{j=1}^{i-1} \wt a_{ij}$ and $c_{i} = \sum_{j=1}^{i}a_{ij}$.
IMEX schemes are a particular case of additive Runge-Kutta methods and so the order conditions can be derived as a generalization of the
notion of Butcher tree; we refer to~\cite{HW} for more details on the order conditions. Before stating the main properties concerning
asymptotic preservation and asymptotic accuracy, we characterize the different IMEX schemes accordingly to the structure of the DIRK method.
Following~\cite{BPR12}, we call an IMEX-RK method of type A if the matrix $A \in \RR^{\nu \times \nu}$ is invertible, or equivalently
$a_{ii}\neq 0$, $i=1,\ldots,\nu$ while we call it of type CK (see \cite{CK}) if the matrix $A$ can be written as
\be
A =
\left(
\begin{array}{ll}
  0  &  0 \\
  a  &  \wh{A}
\end{array}
\right),
\label{CK1}
\ee
with $a=(a_{21},\ldots,a_{\nu 1})^T\in \RR^{(\nu-1)}$ and the submatrix $\wh{A} \in \RR^{(\nu-1)\ \times \ (\nu-1)}$
invertible, or equivalently $a_{ii}\neq 0$, $i=2,\ldots,\nu$.
We write also the matrix $\widetilde A$ for the explicit Runge-Kutta method
\be
\widetilde A=
  \left(
          \begin{array}{ll}
            0 & 0  \\
            \widetilde a & \widehat{\widetilde A}\end{array}
        \right),
\label{CK2}
\ee
where $\widetilde a=(\widetilde a_{21},\ldots,\widetilde a_{\nu 1})^T\in\RR^{\nu-1}$ and $\widehat{\widetilde A}\in\RR^{\nu-1\times\nu-1}$.
We now introduce two useful definitions to characterize the properties of the methods in the sequel.
An IMEX-RK scheme is called \emph{implicitly stiffly accurate (ISA)} if the corresponding DIRK method is \emph{stiffly accurate}, namely
$a_{\nu i}=w_i,\quad i=1,\ldots,\nu$.
If in addition the explicit matrix satisfies $\widetilde a_{\nu i}=\widetilde w_i,\quad i=1,\ldots,\nu$,
the IMEX scheme is said to be \emph{globally stiffly accurate (GSA)} or simply \emph{stiffly accurate}.

Note that for GSA schemes the numerical solution is the same as the last stage value, namely $r^{n+1}=R^{(\nu)}$ and $j^{n+1}=J^{(\nu)}$.
We recall that the order conditions for GSA type $A$ IMEX schemes are particularly restrictive since $\widetilde c\neq c$ and
$\widetilde w\neq w$.
Another restrictive condition is, for high order methods, the request that the matrix $A$ is invertible for details see \cite{DP12, BPR12}.
All these requests make very difficult to derive IMEX GSA schemes of order higher than two.

The detailed analysis of the AP properties of the proposed IMEX schemes is reported in the Appendix. In that part we will give sufficient conditions which guarantee that the schemes are AP and AA. Here, we only recall the main results. Type $A$ IMEX schemes are Asymptotic Preserving and Asymptotically Accurate. If in addition they are GSA the distribution function is projected over the equilibrium at each time step. Two sufficient conditions for type $CK$ IMEX schemes which guarantee the AP and AA properties is that they are $GSA$ and that the initial data are close to the equilibrium state (we say in this case that the initial data are consistent with the limit problem). Again in this case, we get also sufficient conditions to assure that the distribution function is projected over the equilibrium state at each time step.

Here we show how the density values are obtained and that in the limit the schemes solve the drift-diffusion equation with an implicit treatment of the diffusion term. Observe in fact that in order to use schemes (\ref{stRK})-(\ref{numsolj}), we need to know the values of the density distribution $\r$ and of its stages implicitly. These values are obtained by integrating with respect to velocity equation (\ref{stRK}) and (\ref{numsolr}).
We denote by $R=(R^{(k)})_k$, $J=(J^{(k)})_k$ and ${\bf e}=(1,\ldots,1)$ for $k=1,\ldots,\nu$, the column vectors
of the stages for $r$ and $j$ respectively and by $P=(P^{(k)})$ the vector of the stages of the mass density $\r$; it holds
that $P=\int R \,d\bv$.
Moreover, we denote by ${\bf f}_1(R,J) = \Big( f_1( R^{(k)},J^{(k)} ) \Big)_k $ the vector containing stage values for
$k=1,\ldots,\nu$ and similarly for ${\bf f}_2$. The equation for the internal stages in vector form reads
\begin{equation*}
  \int R d\bv = \r^{n} {\bf e} + \Delta t \widetilde{A} \int {\bf f}_1\left(R, J\right)\,d\bv +
           \Delta t A \int {\bf f}_2\left(R\right)\,d\bv,
\end{equation*}
which we can also rewrite as
\begin{equation}\label{Pkkinsvolta}
 P = \r^{n} {\bf e}
   + \Delta t \widetilde{A} \left( \int - \bv\cdot\nabla_\bx \left( J + {\mu\over\l}\bv\cdot\nabla_\bx R \right) \,d\bv
   + \wt G {\bf e} \right) + \Delta t A \mu \Delta_{\bx\bx} \int \frac{|\bv|}\l^2 R \,d\bv,
\end{equation}
while the equation for the numerical solution is
\begin{equation*}
 \r^{n+1} = \r^n + \Delta t \wt w^T \int -\bv\cdot\nabla_\bx \left( J + {\mu\over\l}\bv\cdot\nabla_\bx R \right) \,d\bv
          + \Delta t \wt G + \Delta t w^T \int{\bf f}_2(R) \,d\bv.
\end{equation*}
In the above equation the integral of $R$ is implicit. However, it can be (explicitly) solved by inverting the
matrix describing the discretized diffusion operator related to the term $\int (|\bv|^2/\l) R\,d\bv$. This permits to know the density $\r$ implicitly by the explicit knowledge of $r$ and $j$. Finally, in the limit regime, as shown in the Appendix,
the term $\int (|\bv|^2/\l) R\,d\bv$ is reduced to $D\,P$ which means that we get as desired an implicit discretization of the diffusion term.

\subsection{A linearization technique for the implicit collision term}
In the numerical method described in the previous paragraph the collision operator has to be implicitly computed.
Then, it is necessary that one is able to invert it.
This is usually not the case since, in general, collisions are represented by nonlinear
multidimensional operators which could be costly to compute or even more to invert.
For this reason, we choose to penalize $Q$ with a suitable operator $L$ which $i)$ needs
to be easier to invert and $ii)$ would not change the asymptotic behavior of the solution. The second condition is mathematically expressed by Ker$(Q)=$ Ker$(L)=$ Span$\{M\}$. This strategy has been proposed for the Boltzmann equation in \cite{FilbJin} and subsequently studied in the context of IMEX schemes in \cite{DP12}.

In order to write the modified IMEX schemes, we add and subtract to the collision term $Q$ an operator $L$ and then we combine
the implicit and the explicit solvers as follows:
$$
  \underbrace{Q(r)}_{Implicit}\hspace{4mm}\rightarrow\hspace{4mm}
  \underbrace{\Big( Q(r) - L(r)\Big)}_{Explicit} +
  \underbrace{L(r)}_{Implicit}.
$$
A possible choice for $L$ is a first order approximation of the original operator,
obtained using an expansion of $Q$ near the equilibrium distribution $\r M$:
$Q(r) \approx \nabla_r Q(\r M)( r - \r M )$.\\
Since it is not always possible or easy to compute analytically $\nabla_r Q(\r M)$, one can choose
to approximate it. A possibility is then
\begin{equation}  \label{eq:penbeta}
 L(r) = \beta \, ( \r M - r ),
\end{equation}
where $\beta$ is an upper bound of $\|\nabla_r Q(\r M)\|$.

Regardless from the choice of $L$, we apply the IMEX strategy to the penalized system in the following way
\begin{eqnarray}
  \del_t r & = & \underbrace{ -\bv \cdot \nabla_\bx \left( j + \mu{\bv\over\l}\cdot\nabla_\bx r\right)
             +   E \cdot \nabla_\bv j + {1\over\e^2} \Big( Q(r) - L(r) \Big) + \wh G(r)}_\mathrm{ Explicit } \nonumber \\
           & + & \underbrace{ {1\over\e^2} L(r) +
                 \bv\cdot\nabla_\bx \left(\mu{\bv\over\l}\cdot\nabla_\bx r\right) }_\mathrm{ Implicit },     \label{eq:penr} \\
 \del_t j & = & \underbrace{ -{1\over\e^2} \Big( \l\,j + \bv\cdot\nabla_\bx r - E\cdot\nabla_\bv r \Big) }_\mathrm{Implicit} \label{eq:penj}.
\end{eqnarray}
Observe that computing operator $L$ with the implicit solver stabilizes also the
non-linear collision operator, without changing the asymptotic behavior of the solution. However, this stabilization is not straightforward, on the contrary in order to stabilize the reformulated system it is necessary that the coefficients of the scheme used for the time integration of the linearized collision operator dominate those used for the time integration of the original operator.
Such technique allows us to treat very general collision operators. We prove in the appendix that both $A$ and $CK$ type IMEX schemes if also Globally Stiffly Accurate are AP and AA. More in details, $CK$ schemes needs the additional hypothesis of consistent initial data to assure that the asymptotic properties are satisfied. The request that the schemes are $GSA$ in this case becomes necessary for the stability in the limit of zero mean free path.

\section{Phase-space discretization}
We discuss in this section the discretizations of velocity and space variables. Concerning the velocity variable, because of the particular structure of the problem, it is convenient to
approximate $r$ and $j$ using a Gauss-Hermite expansion.
We decompose then the unknowns $r$ and $j$ as follows \cite{Klar,JP}: $r=\phi M$ and $j=\psi M$, with
$\phi = \phi(t,\bx,\bv)$ and $\psi = \psi(t,\bx,\bv)$.
In this way it is possible to exploit the Gauss-Hermite approximations to efficiently and accurately compute
the derivatives in $\bv$ and the collision operator (which is an integral in $\bv$).
From (\ref{IEsys}) we have:
 \begin{eqnarray}\label{IEphi}
  \del_t(\phi M) & + & \bv\cdot\nabla_\bx \left( (\psi M) + \mu{\bv\over\l}\cdot\nabla_\bx (\phi M)\right)
                   -   E\cdot\nabla_\bv (\psi M) = \nonumber \\
                 & = & {1\over\e^2} Q(\phi M) + \wh G + \mu{|\bv|\over\l}^2\Delta_{\bx\bx} (\phi M),
 \end{eqnarray}
 \begin{equation}\label{IEpsi}
  \del_t(\psi M) = -{ 1\over\e^2} \Big( \l \psi M + \bv\cdot\nabla_\bx(\phi M) - E\cdot\nabla_\bv (\phi M)\Big)\,.
 \end{equation}
From equation (\ref{IEphi}) we obtain
 \begin{eqnarray*}
  M\del_t\phi & + & \bv \cdot \nabla_\bx \left( M\psi + \mu {M\over\l} \bv\cdot\nabla_\bx\phi\right)
                -   E\cdot\left(M\nabla_\bv\psi - {1\over\theta} \bv M \psi \right) = \\
              & = & {1\over\e^2} Q(\phi M) + \wh G + \mu M{|\bv|\over\l}^2\Delta_{\bx\bx}\phi\, ,
 \end{eqnarray*}
from which
 \begin{equation}\label{IEphioM}
  \del_t\phi + \bv\cdot\nabla_\bx \left(\psi + \mu{\bv\over\l}\nabla_\bx\phi\right)
             - E\cdot\Big(\nabla_\bv\psi - {1\over\theta} \bv\psi\Big)
             = {1\over\e^2} \widetilde{Q}(\phi) + G + \mu{|\bv|\over\l}^2\Delta_{\bx\bx}\phi
 \end{equation}
follows, with $ G M = \wh G $ and $ \wt Q M = Q $. Similarly, from equation (\ref{IEpsi}) we have
 \begin{equation}\label{IEpsioM}
  \del_t\psi = - { 1\over\e^2} \left( \l \psi + \bv \cdot \nabla_\bx \phi
               - E \cdot \Big( \nabla_\bv\phi - {1\over\theta} \bv\phi \Big) \right)\,.
 \end{equation}
We conclude this section with an example for $\wt Q$: in the particular case in which $\sigma\equiv1$ in (\ref{eq:Q}),
we get the co called \emph{relaxed time approximation} (RTA). It is easy to see then that
$$
  Q(f)=\r M - f \quad \Rightarrow \quad Q(\phi M) = \rho M - \phi M = \Big(\rho-\phi\Big) M
$$
and thus it holds that $\quad \wt Q(\phi) = \rho-\phi$.

\subsection{Velocity discretization}
We describe the Gauss-Hermite approximation (\cite{Ch, Kl2, JP}) in the monodimensional case.
The multidimensional case is obtained applying the monodimensional rule dimension-by-dimension.

Let consider $r=\phi M$ and $j=\psi M$, with
$$
  \phi(\bv) = \sum_{k=0}^N \phi_k \tH_k(\bv) \,, \qquad
  \psi(\bv) = \sum_{k=0}^N \psi_k \tH_k(\bv) \,,
$$
being the Hermite expansion. Here $\tH_k$ are the renormalized Hermite polynomials and coefficients $\phi_k$ and
$\psi_k$ can be computed thanks to the inverse expansion (we refer to \cite{Kl2} for more details).
The computation of the collision operator becomes
\begin{eqnarray*}\label{qcoll}
  Q(r)(\bv)  & = & M(\bv) \sum_{j=0}^N \sigma(\bv, \bv_j) \, \phi(\bv_j) \, \bw_j - \lambda(\bv) \, r(\bv)\,, \\
  \text{with }\quad \lambda(\bv) & = & \sum_{j=0}^N \sigma(\bv, \bv_j) \, \bw_j\,,
\end{eqnarray*}
where $(\bv_j, \bw_j)$ are points and weights of the Gauss-Hermite quadrature rule.
Finally, the derivatives with respect to $\bv$, which are given by
$$
  \nabla_\bv r = M \nabla_\bv \phi - \frac 1 {\theta} \bv M \phi \,, \quad
  \nabla_\bv j = M \nabla_\bv \psi - \frac 1 {\theta} \bv M \psi \,,
$$
become
\begin{equation}     \label{eq:derv}
 \nabla_\bv \phi =
 \sum_{j=0}^N \phi(\bv_j) c_j(\bv), \quad \text{ and } \quad
 \nabla_\bv \psi = \sum_{j=0}^N \psi(\bv_j) c_j(\bv),
\end{equation}
$$
  c_j(\bv) =  \sum_{k=1}^N \sqrt{2k}\, \tH_k(\bv_j) \, \tH_{k-1}(\bv) \, \bw_j \,.
$$
\smallskip\noindent
{\bf Remark:} Coefficients $c_j(\bv_i)=c_{ij}$ for any component of $\bv$ can be computed at the beginning of the simulation and stored in a matrix since they do not depend on functions $\phi$ and $\psi$.

\subsection{Space discretization}

In this section we emphasize some requirements about the space discretization of the system.
We want our scheme to work both in the kinetic regime ($\e\gg0$), in which the hyperbolic behavior
is more relevant, and in the limit regime ($\e\approx0$), in which the system is characterized by diffusive behavior.
Moreover, the characteristic speeds of the system (which are of the order of $1/\e$) tend to infinity as $\e\rightarrow0$ and so shock
capturing methods based on characteristics directions, such as, e.g., upwind methods, become useless.
On the other hand, central differences schemes avoid excessive dissipation but, when $\e$ is not small or when the limiting equations
contain advection terms, may lead to unstable (or not accurate) discretizations.

In order to overcome these well-known facts and to have the correct asymptotic behavior, we fix
some general requirements for the space discretization:
\begin{description}
 \item[i)] \emph{correct diffusion limit:} as we already observed in previous section, if we want a correct approximation in the limit
  case $\e=0$, we need that $\mu(0)=1$ and to use the  same space discretization for the transport terms in (\ref{IEphioM})
  and (\ref{IEpsioM});
 \item[ii)] \emph{compact stencil:} we want to use a scheme with a compact stencil in the diffusion limit $\e\rightarrow0$.
  This property is satisfied  if point {\bf i)} is satisfied and we use a suitable discretization for the second order derivative that
  characterizes the diffusion limit;
 \item[iii)] \emph{shock capturing:} the chosen scheme should be based on high order shock capturing fluxes for the convection part.
  This is necessary not only for large values of $\e$ but also when we consider convection-diffusion type limit equations with small
  diffusion. The high order fluxes are then necessary for all space derivatives except for the second order term
  $\mu\,|\bv|^2\Delta_{xx} r/\l$ on the right hand side in (\ref{IEphioM});
 \item[iv)] \emph{avoid solving nonlinear algebraic equations:} in order to have a more efficient method, we do not want to solve
  the nonlinear equation which comes from the implicit treatment of the space derivative in equation (\ref{IEpsioM}) for $\psi$,
  (and in (\ref{eq:penj}) for the odd parity $j$). To achieve this we have chosen a partitioned approach for the time integration, thanks to which the values
  of $\phi$ are already available from the previous solution of (\ref{IEphioM}).
\end{description}

\subsubsection{Modified fluxes}

In order to satisfy the above requirements we choose to use a Lax-Friedrichs type flux with high order WENO reconstruction \cite{Shu, CGMS1, CGMS2}.
This gives us the ability to ensure accuracy and also to stabilize the solution in the presence of discontinuities
or arising shocks.

Our strategy is the following: in the kinetic regime, where transport dominates the dynamics, we use the standard
Lax-Friedrichs type flux with WENO reconstruction for the derivatives \cite{BPR12,JPT1,NP1}.
In this case, indeed, this is a proper strategy which allows us to obtain also high order accuracy.
When considering the limiting regime instead, we cannot use the Lax-Friedrichs scheme as it is.
As we will show later in this section, the numerical viscosity introduced by such scheme in this case is proportional
to $1/\e$.
Clearly, when the mean free path goes towards zero ($\e\rightarrow0$) such quantity is too large and
causes loss of accuracy (see for instance \cite{NP1, NP2}).
Thus, in such situation we decide to bound the numerical viscosity modifying the fluxes.
This is possible because when $\e$ becomes small the diffusive regime becomes dominant and thus stability is granted by the ``physical''
viscosity given by the system itself. 
We point out that, in this work, the stability of the proposed modified fluxes approach is supported by numerical evidence.
Theoretical estimates for this technique and for this kind of scaling will be the subject of a future work.

We present here the modified fluxes approach using a model problem. Given the transport equation for the unknown $w=w(t,x)$
$$
  w_t + \del_x a(w) = 0,
$$
with $a$ a hyperbolic flux, we consider a complete discretization in conservation form
$$
  w_i^{n+1} = w_i^n - \lambda \Big( W_{i+1/2} - W_{i-1/2} \Big),
$$
\begin{equation} \label{eq:numflux}
 W_{i+1/2} = {1\over2}\Big[ a(w_{i+1}^n) + a(w_i^n) - \alpha \, ( w_{i+1}^n - w_i^n ) \Big],
\end{equation}
where $ \lambda = \Delta t / \Delta x $ and the parameter $\alpha$ represents the numerical viscosity.
The standard Lax-Friedrichs scheme requires the value $ \alpha = 1 / \lambda $.
Such scheme is stable if the following two inequalities are satisfied
\begin{equation}\label{lambdalpha}
\lambda \leq { 1 \over \max_w |a^\prime(w)| } \quad \text{ and } \quad
\max_w|a^\prime(w)| \leq \alpha \leq { 1 \over \lambda }.
\end{equation}
Formula (\ref{eq:numflux}) is the basis for the construction of a conservative numerical flux which can
be implemented using ENO or WENO high order reconstructions \cite{Shu}.

To derive the modified fluxes, let consider now the prototype system
\begin{equation}\label{eq:1}
\begin{split}
u_t       & =  - ( v - \mu u_x )_x + \mu u_{xx}\\
\e^2 v_t  & =  u - u_x - v
\end{split}
\end{equation}
which shares the same structure of our original problem. In the limit $\e=0$, if $ \mu(\e=0) = 1$, the above system leads to the drift-diffusion equation $ u_t + u_x = u_{xx}$.
We then write a semi-discrete approximation of (\ref{eq:1}) as
$$
(u_i)_t = - \frac{ 1 }{ \Delta x } \Big( U_{i+1/2} - U_{i-1/2} \Big),
\qquad
\e^2 (v_i)_t = - \frac{ 1 }{ \Delta x } \Big( V_{i+1/2} - V_{i-1/2} \Big) + u_i - v_i,
$$
with numerical fluxes given by
$$
U_{i+1/2} = {1\over2} \left[ \Big( v_{i+1} + v_i \Big ) - \alpha \, \Big( u_{i+1} - u_i \Big ) \right],
\quad
V_{i+1/2} = {1\over2} \left[ \Big( u_{i+1} + u_i \Big ) - \e^2 \alpha \, \Big( v_{i+1} - v_i \Big ) \right]
$$
where $\alpha$ is the numerical viscosity as before. Rewriting the numerical fluxes as
\begin{equation} \label{eq:modflumix}
U_{i+1/2} = {1\over2} \left[ \Big( v_{i+1} + v_i \Big ) - \alpha_u \, \Big( u_{i+1} - u_i \Big ) \right],
\quad
V_{i+1/2} = {1\over2} \left[ \Big( u_{i+1} + u_i \Big ) - \alpha_v \, \Big( v_{i+1} - v_i \Big ) \right]
\end{equation}
we have
\begin{equation}\label{eq:fluxphys}
    \alpha_u = \frac{1}{\e},
    \qquad \qquad
    \alpha_v = \e.
\end{equation}
The stability conditions (\ref{lambdalpha}) in this case read:
$$
  \lambda = { \Delta t \over \Delta x } \leq \e \hspace{5mm} \Rightarrow \hspace{5mm} \Delta t \leq \e \, \Delta x,
$$
$$
\frac{ 1 }{ \e } \leq \alpha \leq \frac{ 1 }{ \lambda } \hspace{5mm} \Rightarrow \hspace{5mm} \alpha \geq \frac{ 1 }{ \e }.
$$
As pointed out before, we want to avoid such restrictive time step. 
To this aim, close to the limit we modify the numerical viscosity setting
\begin{equation}\label{eq:fluxnum}
    \alpha_u = 1,
    \qquad  \qquad
    \alpha_v = \e^2.
\end{equation}
More in details, until we reach a regime in which the physical diffusion is not large enough to guarantee stability
we need to satisfy~(\ref{eq:fluxphys}).
On the other hand, when physical diffusion becomes relevant, we can avoid such restrictions and we can choose the modified
fluxes~(\ref{eq:fluxnum}). The practical choice we did in our numerical tests is
\begin{equation}\label{eq:fluxnum_used}
    \alpha_u = \min(\frac{1}{\e},1),
    \qquad  \qquad
    \alpha_v = \min(\e,\e^2).
\end{equation}

To summarize, we write here the complete numerical discretization of system (\ref{eq:penr})-(\ref{eq:penj}).
We denote the values $\phi_{ij}^n = \phi(t^n,\bx_i,\bv_j)$ for $i$ and $j$ varying in the phase-space index set and at time $t=t^n$.
We denote in the same way the other variables appearing in the sequel. For stage vectors $\Phi$ and $\Psi$ we denote
$\Phi_{ij} = ( \Phi_{ij}^{(k)} )_k $ and $\Psi_{ij} = ( \Psi_{ij}^{(k)} )_k $ for $k=1,\ldots,\nu$. Then we have
\begin{eqnarray}
 \Phi_{ij} & = & \phi_{ij}^n {\bf e} + \Delta t \, \wt A \Big[ - \Gamma_{ij}( \Psi^*, \Phi, \alpha_v )
               + E_i \otimes ( \Psi_{ij}^\bv - 2\bv_j\Psi_{ij} ) + G_{ij} {\bf e}                            \nonumber  \\
           & + &  \frac 1 {\e^2} \Big( {\bf {\wt Q}}_{ij}( \Phi ) - {\bf {\wt L}}_{ij}(\Phi) \Big) \Big]
               + \Delta t \, A \left[ {1\over\e^2} {\bf {\wt L}}_{ij}( \Phi )
               + \frac{\mu}{\l_j} \, \bv_j^2 \, \Phi_{ij}^{\bx\bx} \right],                                 \label{eq:dfPhhi}  \\
 \Psi_{ij} & = & \psi_{ij}^n {\bf e} - \frac{\Delta t}{\e^2} \, A \Big[ \l_j \Psi_{ij}
          + \Gamma_{ij}( \Phi, \Psi, \alpha_u ) - E_i \otimes ( \Phi_{ij}^\bv - 2 \bv_j \Phi_{ij} ) \Big]. \label{eq:dfPssi}
\end{eqnarray}
and $\phi_{ij}^{n+1} = \Phi_{ij}^\nu$ and $\psi_{ij}^{n+1} = \Psi_{ij}^\nu$ for the numerical solution.
To make formulas more readable we defined some shorthands.
Operator $\Gamma$ stands for the numerical discretization of the transport derivatives: $\Gamma\approx\bv\cdot\nabla_\bx$,
obtained as in (\ref{eq:modflumix}) and it reads
$$
\Gamma_{ij}(h,k,\alpha) = \frac{1}{\Delta x} \Big( H_{i+\frac12} - H_{i-\frac12} \Big),
$$
with
\begin{equation}
 H_{i+\frac12} = \frac{\bv_j}{2}\Big( h_{i+1} + h_i - \alpha ( k_{i+1} - k_i ) \Big).
 \label{eq:hf}
\end{equation}
The term $\Psi^*$ is $\Psi^* = \Psi + ( \mu / \l ) \Gamma( \Phi, \Psi, \alpha_u )$ while $\Phi_{ij}^\bv$ and $\Psi_{ij}^\bv$
stand for the derivative with respect to $\bv$ and are obtained by means of (\ref{eq:derv}).

The diffusion term in the r.h.s. of equation (\ref{eq:dfPhhi}), i.e. the derivative $\Phi^{\bx\bx}$, stands for the standard
central second order finite difference technique, i.e.
$$
  \Phi_{ij}^{\bx\bx} =  \frac 1 { \Delta x^2 } \Big( \Phi_{i+1,j} - 2\Phi_{ij} + \Phi_{i-1,j} \Big),
$$
when the second order time discretizations is used, and for the standard central fourth order finite difference technique when the
third order time discretizations are used.
Regarding the numerical viscosity in equations (\ref{eq:dfPhhi}) and (\ref{eq:dfPssi}), if we are in the kinetic regime
we choose the physical values given by (\ref{eq:fluxphys}) while in the limiting regime we consider the modified ones
given by (\ref{eq:fluxnum}).
The space discretization of the electric potential $\varPhi$ can be performed by standard methods \cite{Kl2, SZ}.

\subsection{Boundary conditions}
The treatment of boundary conditions for the Boltzmann-Poisson problem in the diffusive limit is, in the general case, a very hard task.
It is necessary to tackle several difficulties, such as complex geometries of the boundaries and to take into account boundary layers.
For a consistent treatment of boundary condition see, for instance, \cite{Kl2,LemMeh2012,LemMieu2008} and references therein.
All these aspects are out of the scope of this article and we will only deal with assigned, constant in time boundary data,
the so called maxwellian injection.
We show now a possible strategy for a consistent numerical implementation of such conditions in the one-dimensional situation.\\
For $x \in\, (x_L, x_R)$ a max\-wel\-lian injection is defined by
$$
  f(t,x_L,v) = F_L ( v ),\qquad
  f(t,x_R,-v) = F_R ( v ),
$$
for $v>0$, where $F_L$ and $F_R$ are assigned nonnegative functions proportional to the max\-wel\-lian distribution $M$.
We numerically approximate these conditions in two different ways, depending on the regime in which the system is.

In the kinetic regime, for $x=x_L$ we set $f(t,x_L,v) = F_L(v)$, we extrapolate $f(t,x_L,-v)$ (the outgoing particles)
from the values of $r$ and $j$ inside the domain and then we define $r(x_L)$ and $j(x_L)$ thanks to the parity formulas
(\ref{eq:parir}) and (\ref{eq:parij}). At the right boundary $x_R$ a similar treatment is used.

In the diffusive regime instead, to get a boundary condition for $r$ and $j$ we use the relations (for positive $v$ only)
\begin{equation}
 r + \e j \, \big|_{x=x_L} = f(t,x_L,v) = F_L\,, \qquad
 r - \e j \, \big|_{x=x_R} = f(t,x_L,-v) = F_R\,.
 \label{thbound}
\end{equation}
Then we consider equation (\ref{thq:2}), which gives a good approximation of $j$ when $\e$ is small,
i.e. $\lambda j = - \bv \nabla_\bx r + E \nabla_\bv r$, and applying it in (\ref{thbound}) one gets
$$
  r- \frac{\e}{\lambda}(\bv \nabla_\bx r - E \nabla_\bv r) \big|_{\bx=\bx_L} = F_L\,, \qquad
  r+ \frac{\e}{\lambda}(\bv \nabla_\bx r - E \nabla_\bv r) \big|_{\bx=\bx_R} = F_R\,.
$$
To approximate $\nabla_\bx r(\bx_L)$ and $\nabla_\bx r(\bx_R)$ we use one-sided finite difference discretizations
of the desired order and to approximate $\nabla_\bv r$ we observe that from (\ref{thbound}) it holds that
$$
\nabla_\bv r \big|_{\bx=\bx_L} = \nabla_\bv F_L + O(\e), \qquad
\nabla_\bv r \big|_{\bx=\bx_R} = \nabla_\bv F_R + O(\e),
$$
which in the end leads to (up to $O(\e^2)$)
\begin{equation}
 \label{bcok}
  r- \frac{\e}{\lambda}(\bv \nabla_\bx r - E \nabla_\bv F_L)\big|_{\bx=\bx_L} = F_L\,, \qquad
  r+ \frac{\e}{\lambda}(\bv \nabla_\bx r - E \nabla_\bv F_R)\big|_{\bx=\bx_R} = F_R\,.
\end{equation}

\section{Numerical tests}
In this section, we present several numerical results to test the performance of the proposed schemes.
We show that our schemes are computationally very efficient, i.e. $\Delta t = O(\Delta x)$ while at the same time they preserve
high order of accuracy in all regimes.
The setting is a monodimensional phase-space, i.e. $x,v\in\RR$.

In our computations we use two different scattering cross-sections:
a simple isotropic case with a constant cross-section $\sigma^{RTA}(\bv,\bw)=1$, this corresponds to the relaxation time
approximation (RTA), for which the collision operator has the simple form
$$
  Q(f)=\rho M-f
$$
and a regularized anisotropic cross section \cite{MPS} for electron-phonon interactions (EPI)
$$
  \sigma^{EPI}(\bv,\bw)=\wt{\delta}(|\bv|^2-|\bw|^2+1) + \wt{\delta}(|\bv|^2-|\bw|^2-1),
$$
where $\wt\delta(x)=\exp(-C|x|^2)$ is a smoothed delta function with $C$ a positive constant (we set $C=1/10$).
When we consider the EPI model, we apply the penalization technique using operator $L$ given by (\ref{eq:penbeta}) with $\beta=1$ (observe that this corresponds to nothing else but the RTA approximation).

Concerning the value of $\mu$, we choose a simple form given by
\begin{equation}
 \mu(\e,\Delta x) =
 \left\{
 \begin{array}{lll}
  1, & \mbox{if} & \e<\Delta x,\\
  0, & \mbox{if} & \e\geq\Delta x.
 \end{array}
 \right.
\end{equation}
As we already observed, in this work we assume $\mu$ has a constant for a given value of $\e$.
We recall here that more accurate choices are possible for $\mu$, i.e. $\mu=\mu(\e,\Delta x, \Delta t)$, and we remind to a future work
for a deeper analysis of this aspect.

In all our numerical tests the discretization in the velocity space is obtained using $N_v$ Gauss-Hermite quadrature points,
with $N_v=16$: there are 8 nodes for positive velocities and 8 for negative ones (we used scaled values in order to consider
the range $[-v_{\max},v_{\max}]$, with $v_{\max}\approx5$).
The influence of the number of quadrature points on the accuracy of the results is treated, for instance, in \cite{Kl2}. The IMEX schemes we used for our simulations are the second order IMEX ARS-(2,2,2) scheme \cite{ARS} and
the third order IMEX BPR-(3,5,3)  scheme \cite{BPR12}. For the sake of completeness we report the Butcher tables of the schemes in the Appendix. We compare our results also with a simple first order IMEX scheme, obtained by combining the first order implicit and explicit Euler schemes. A reference solution is always reported for all the tests.

\subsection*{Test \# 1}
In this problem we have a potential well in the left half of the slab. We test the behavior of the scheme when the system
is subject to a constant in time, electric field which varies along the $x$-axis: it ranges from a minimum value of -10
to a maximum value of $10$.
We perform simulations in both kinetic and fluid regimes, using both scattering kernels. \\
In the kinetic regime, i.e. $\e=1$, we stop our simulations at time $T_f=0.08$, with $x\in [0,1]$ using $50$ grid points and with
an initial distribution given by $f(x,v,t=0)=M(v)$.
At the boundaries we set the values $F_L(v) = M(v),$ and $F_R(v)=M(v)$
and we approximate them as described in previous section. The other parameters of the simulations are
$$
G=0, \quad \vpp = \exp(-c(1/4-x)^2),
$$
with $c=50\exp(1)$. Since we are in the kinetic regime, the time step is given by the hyperbolic condition
\begin{equation}  \label{eq:dtH}
 \Delta t = \Delta t_H = c_H \, \e \, \Delta x/v_{\max}.
\end{equation}
For this test we set the CFL constant to $c_H = 0.5$ for all schemes.

\begin{figure}[ht!]
 \centering
 \includegraphics[bb=0 0 360 252,scale=1.0,keepaspectratio=true]{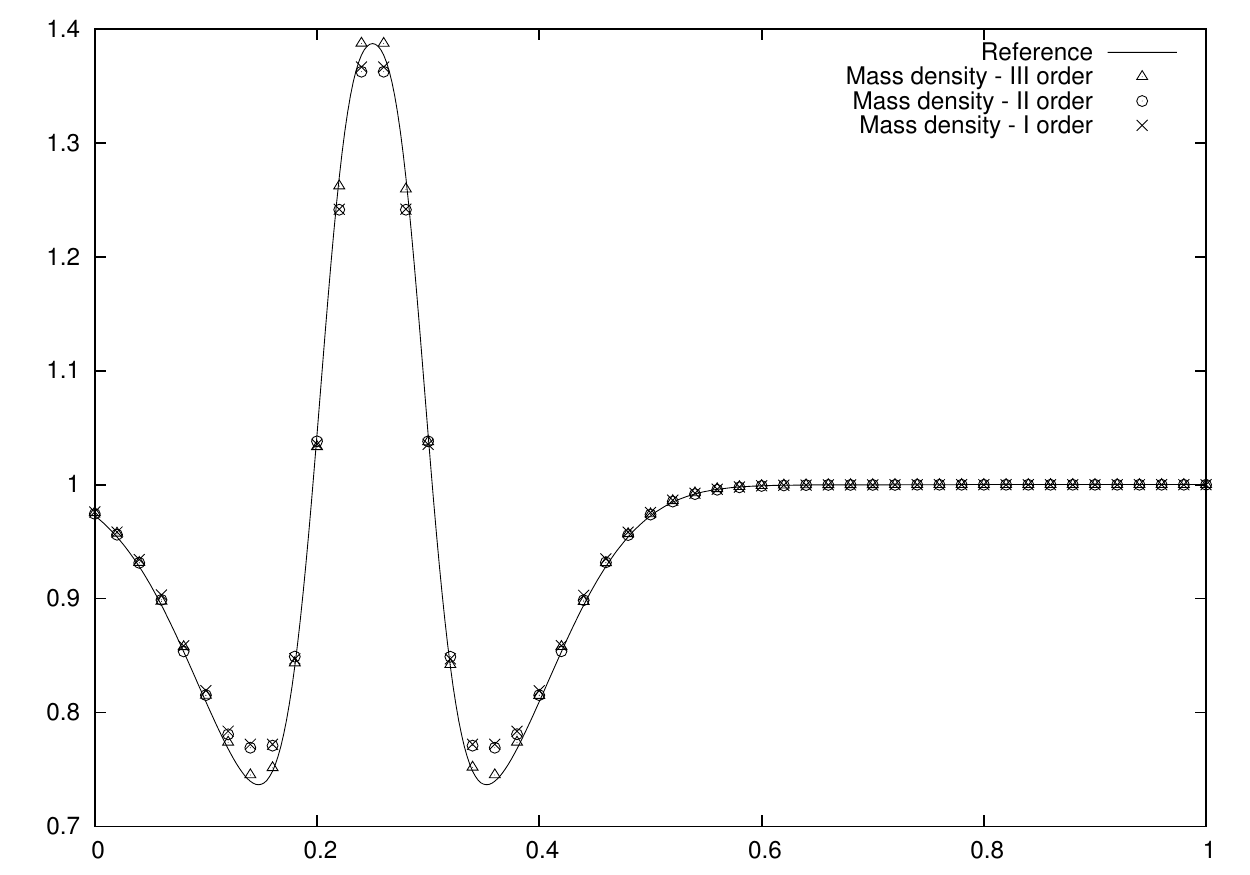}
 \caption{Test \# 1: comparison of the reference (line) and numerical ( ($\times$) I order, ($\circ$) II order and ($\triangle$) III order)
  mass distributions in the kinetic regime  $\e=1$ with potential well for the RTA model, $\Delta t = \Delta t_H$.
  On the $x$-axis the space variable, on the $y$-axis the mass $\r$.}
 \label{fig:test0_RTA}
\end{figure}

The computed solutions of the mass density are presented in figure \ref{fig:test0_RTA}. We show results for the RTA kernel since
the obtained results for the EPI model are similar.
We compare the reference density, obtained with a fourth order explicit RK scheme with third order WENO reconstruction using $N_x=400$,
with the first, second and third order IMEX approximations. As expected, the third order scheme gives a more accurate solution.

\begin{figure}[ht!]
 \centering
 \includegraphics[bb=0 0 360 252,scale=1,keepaspectratio=true]{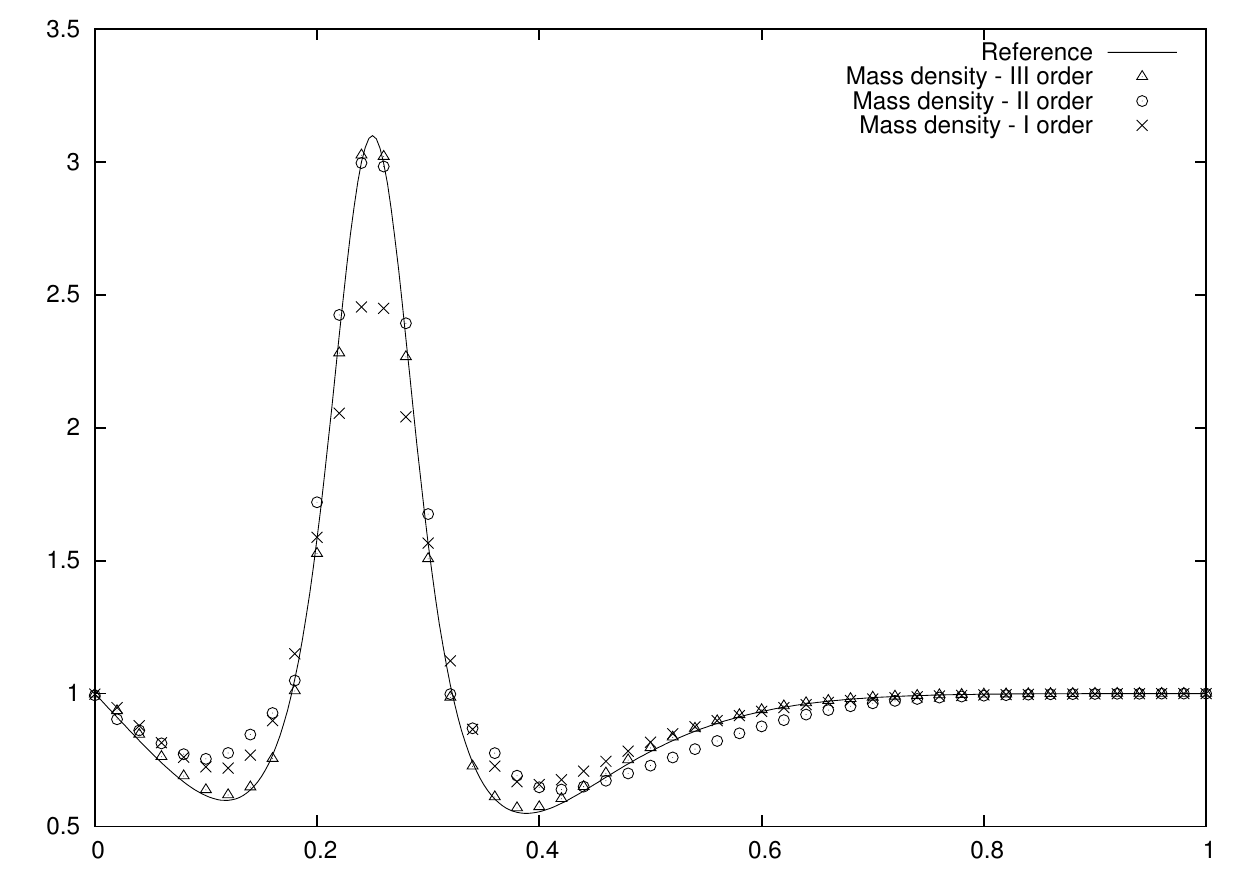}
 \caption{Test \# 1: comparison of the reference (line) and numerical ( ($\times$) I order, ($\circ$) II order and ($\triangle$) III order)
  mass distributions in the limit regime  $\e=0.002$ with potential well for the RTA model, $\Delta t = \Delta t_M$.
  On the $x$-axis the space variable, on the $y$-axis the mass $\r$.}
 \label{fig:test1_RTA}
\end{figure}

Next in figure \ref{fig:test1_RTA} we report the results obtained in the fluid regime in the RTA case.
We set in this case $\e=0.002$ and stop the simulation at $T_f=0.03$ using $50$ grid points while the other parameters are the same as in
the kinetic test case. Now the time step is
\begin{equation} \label{eq:dtm}
 \Delta t = \Delta t_M = c_M \Delta x,
\end{equation}
with $c_M=0.5$ for all orders. The reference solution in this case is obtained as in \cite{JP} with $N_x=400$.

Results clearly show that we are able to overcome the very severe time step restriction that would be necessary for
this kind of test, using standard techniques, while still attaining high order of accuracy.
Indeed, a simple explicit integrator would require
$$
  \Delta t = \min\left\{ \Delta t_P = \frac{\Delta x^2}{2}, \, \Delta t_H  = c_H \, \e \, \Delta x/v_{\max}\right\}
$$
with $\e=2\cdot10^{-3}$. On the other hand, an IMEX strategy without the implicit treatment of the diffusive term would require
$\Delta t = \Delta t_P$, which is $\Delta x$ times smaller then ours.
Indeed, we are able to compute the solution with only a linear dependence of $\Delta t$ on $\Delta x$ which is given by (\ref{eq:dtm}).

\subsection*{Test \# 2}
\begin{figure}[ht!]
 \centering
 \includegraphics[bb=0 0 360 252,scale=1,keepaspectratio=true]{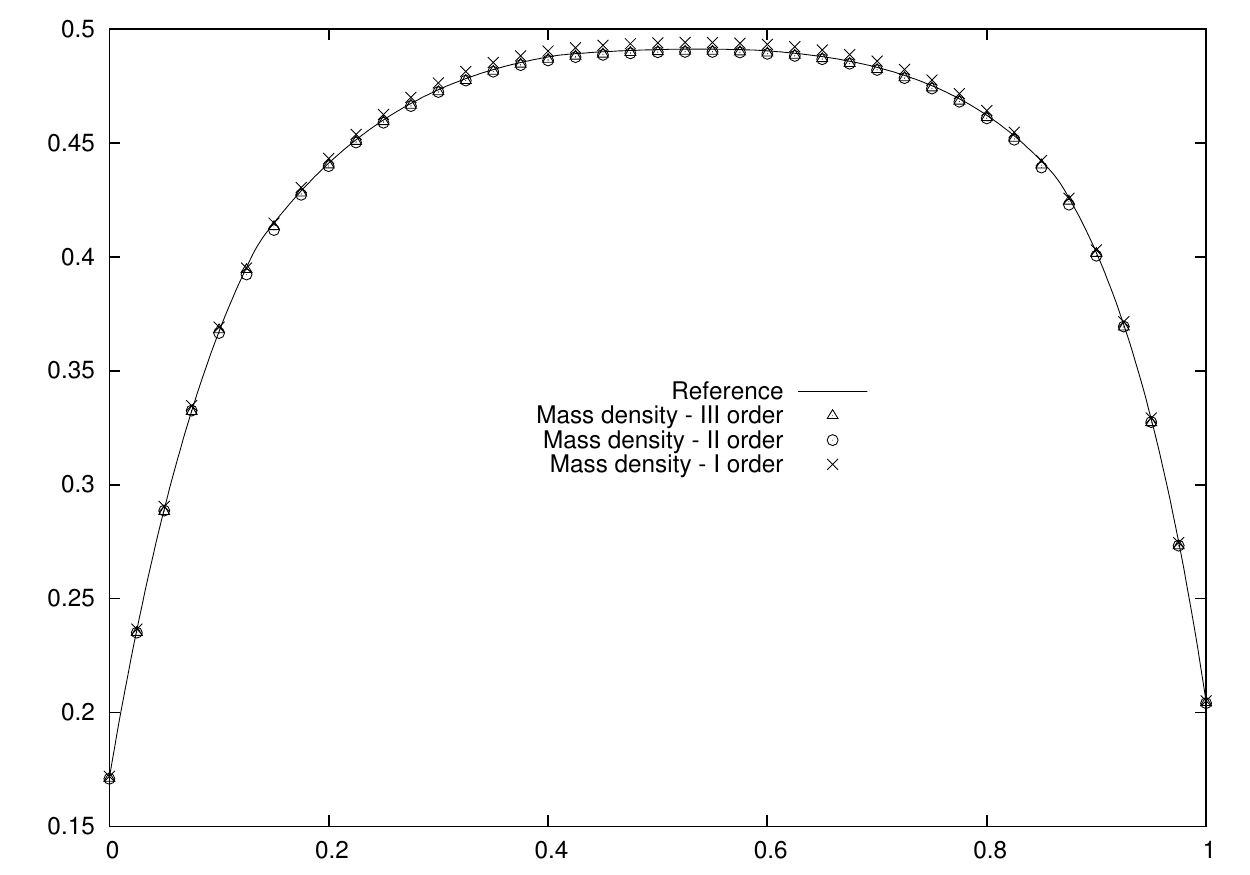}
 \caption{Test \# 2: comparison of the reference (line) and numerical ( ($\times$) I order, ($\circ$) II order and ($\triangle$) III order)
  mass distributions in the kinetic regime $\e=1$ with a constant electric field and non zero source term for the EPI model, $\Delta t = \Delta t_H$.
  On the $x$-axis the space variable, on the $y$-axis the mass $\r$.}
 \label{fig:test2_EPI}
\end{figure}

\begin{figure}[ht!]
 \centering
 \includegraphics[bb=0 0 360 252,scale=1,keepaspectratio=true]{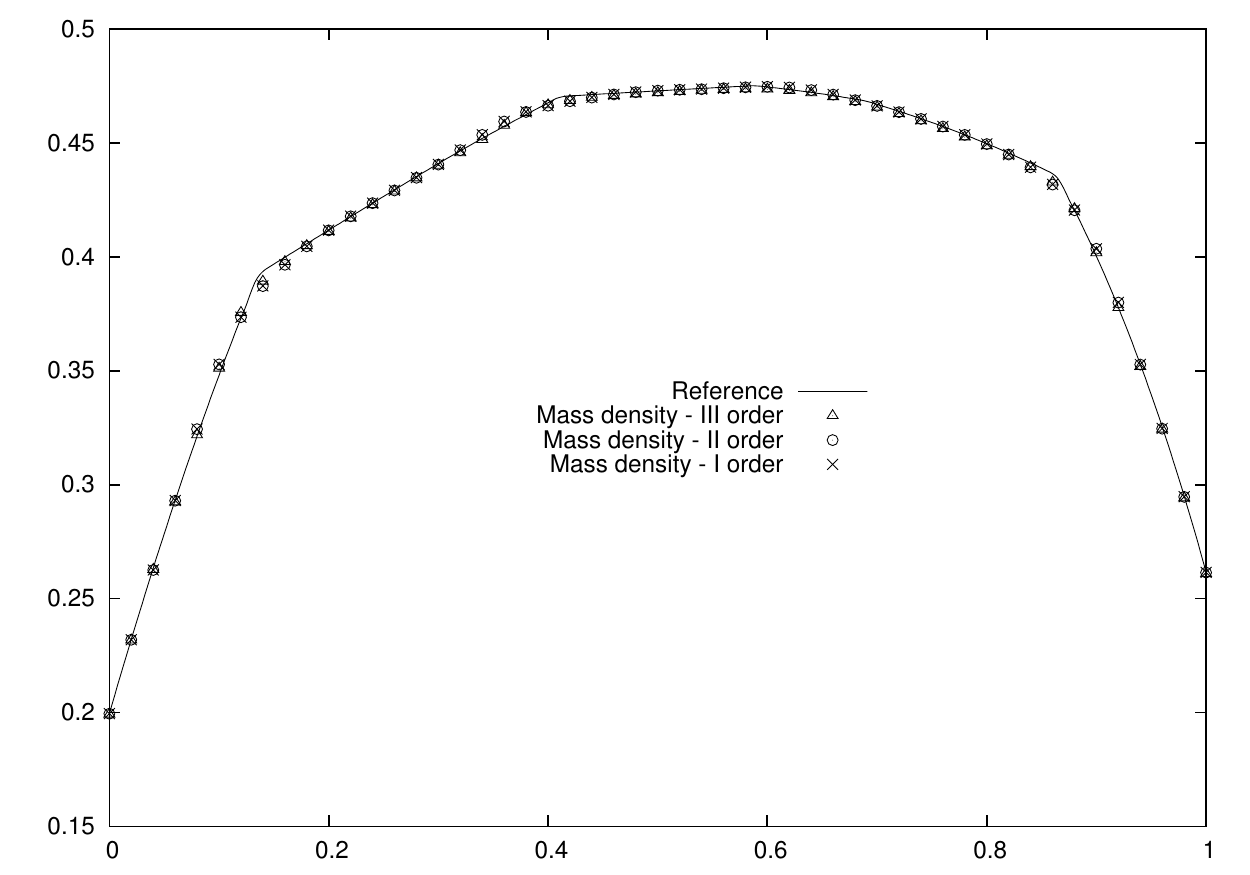}
 \caption{Test \# 2: comparison of the reference (line) and numerical ( ($\times$) I order, ($\circ$) II order and ($\triangle$) III order)
  mass distributions in the kinetic regime $\e=1$ with a constant electric field and non zero source term for the RTA model, $\Delta t = \Delta t_H$.
  On the $x$-axis the space variable, on the $y$-axis the mass $\r$.}
 \label{fig:test2_RTA}
\end{figure}

In this test we study the behavior of the scheme when a non zero source term is present.
The system consists of a diffusive slab with a flat interior source and a constant electric field $E=-1$.
The setting of the problem is as follows: $x\in [0,1]$, $\, \vpp=x$ and $\, G=1$.
The initial distribution is $f(x,v,t=0)=0$ and at the boundaries we set $F_L(v) = 0$ and $F_R(v)=0$.
We perform the simulation in the kinetic regime, i.e. $\e=1$, and we stop computations at time $T_f=0.5$, using $50$ grid points.
The solution for the EPI model is given in figure \ref{fig:test2_EPI} while that for the RTA model is given in figure \ref{fig:test2_RTA}.
For both tests the hyperbolic time step condition $\Delta t_H = c_H \, \e \, \Delta x/v_{\max}$ is imposed.
The reference solution is again obtained with a fourth order explicit RK scheme and WENO reconstruction with $N_x=400$.
The Figures show the good behavior of the numerical schemes. In the EPI case it is possible to see the difference of the first order
from the higher order schemes.
A closer look at the results reveals that higher order methods guarantee more precision with a similar time step condition also
in the RTA case.

\begin{figure}[ht!]
 \centering
 \includegraphics[bb=0 0 360 252,scale=1,keepaspectratio=true]{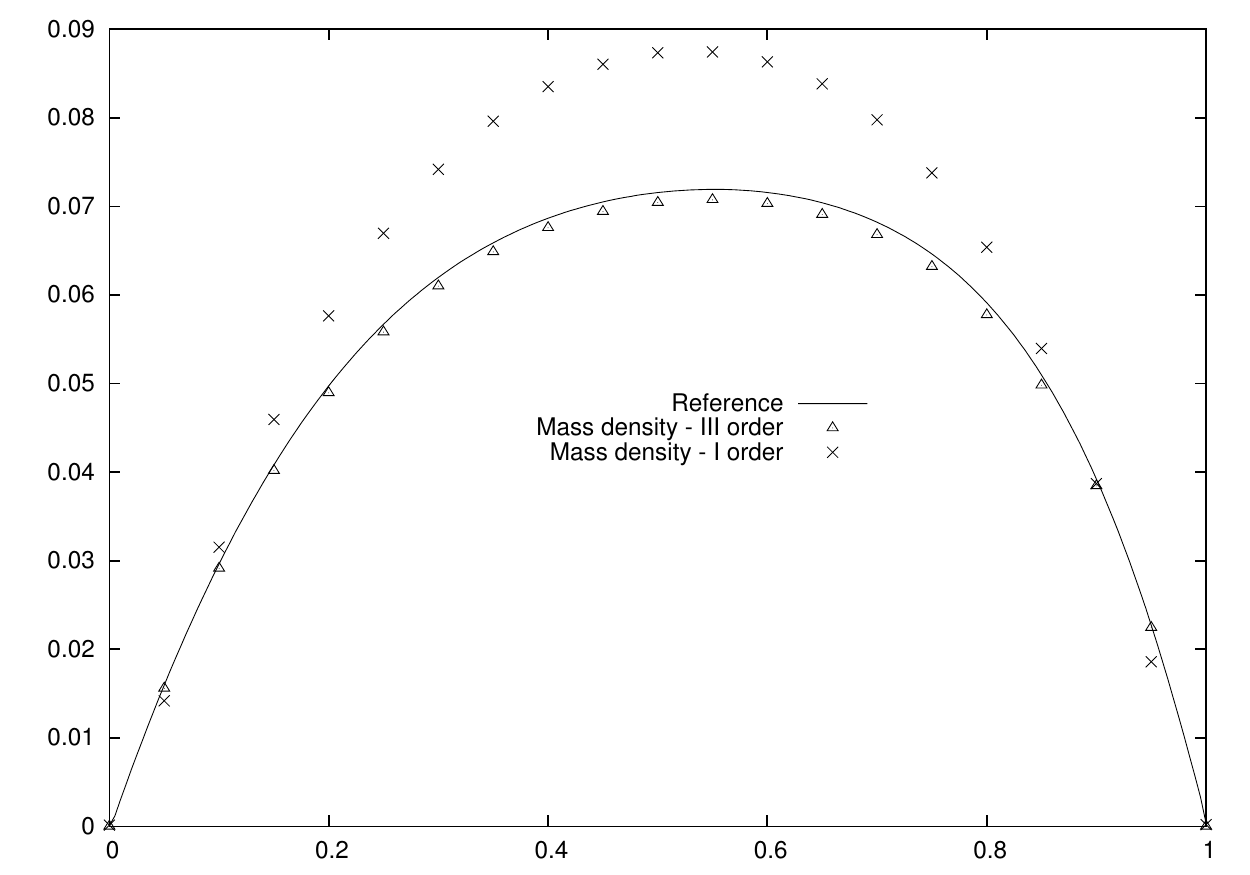}
 \caption{Test \# 2: comparison of the reference (line) and numerical ( ($\times$) I order and ($\triangle$) III order)
  mass distributions in the limit regime $\e=0.001$ with a constant electric field and non zero source term for the RTA model, $\Delta t = \Delta t_M$.
  On the $x$-axis the space variable, on the $y$-axis the mass $\r$.}
 \label{fig:test3_RTA}
\end{figure}

Concerning the limit regime, i.e. $\e=0.001$, the solution at time $T_f=0.1$ for the RTA model, using $20$ grid points and a time step
$\Delta t = c_M \Delta x$, again with $c_M=0.5$ for the three methods is reported in figure \ref{fig:test3_RTA}.
The reference solution in this case is obtained as in \cite{JP} with $N_x=200$, which is $10$ times finer then our grid.
One can clearly observe that even for coarse discretizations the behavior of the system is well described and that the third order scheme gives very accurate results. We point out that in this problem the use of an high order in time method gives results which are much more closer to the reference solution than those obtained by a first order in time scheme. This confirms the importance of developing high order schemes for multiscale problems.

\subsection*{Test \# 3}
In this last test, we consider a unipolar diode of type $\rho^+ \rho \rho^+$ in the
diffusive regime \cite{JP}.
In this case the electric field is self-consistently computed by the solution of the Poisson equation :
\begin{equation} \label{firo}
\gamma \Delta_x \vpp = \rho -\rho_d(x), \qquad \vpp(0)=0,\quad \vpp(1)=V,
\end{equation}
where $\gamma$ is the scaled Debye length, $V$ is the voltage applied at the right boundary and $\r_d(x)$ is the doping profile
$$
\rho_d(x) = 1 - \frac{1-m}{2} \left[\tanh\left(\frac{x-x_1}{s}\right) - \tanh\left(\frac{x-x_2}{s}\right)\right],
$$
with $s=0.02$ (which controls the thickness), $m=0.001$ (the minimum value), $x_1=0.3$ and $x_2=0.7$.
We set in addition $\e=0.001$, $f(x,v,t=0)=M(v)$, $N_x=50$, while the time step is $\Delta t = c_M \Delta x$,
with $c_M = 0.1$. The other parameters are
$$
x\in [0,1], \quad G=0, \quad F_L(v) = M(v), \quad F_R(v)=M(v).
$$
According to \cite{FJO} at the boundary we assume
$$
\partial_x j(0,t) = \partial_x j(1,t)=0.
$$
\begin{figure}[ht!]
 \centering
 \includegraphics[bb=0 0 360 252,scale=1,keepaspectratio=true]{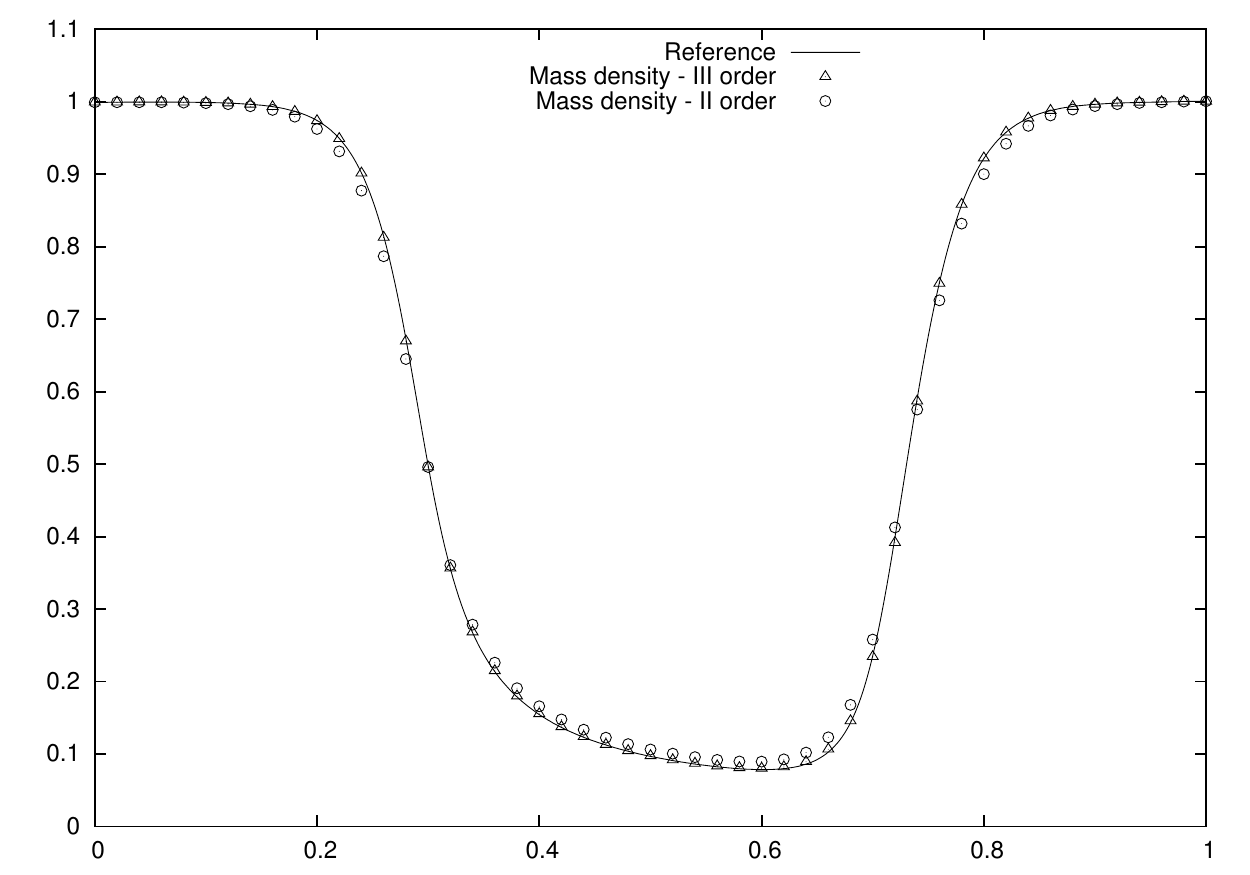}
 \caption{Test \# 3: comparison of stationary reference (line) and numerical ( ($\times$) I order, ($\circ$) II order and ($\triangle$) III order)
  mass distributions for the Vlasov-Poisson problem in the limit regime $\e=0.001$ for the RTA model, $\Delta t = \Delta t_M$.
  On the $x$-axis the space variable, on the $y$-axis the mass $\r$.}
 \label{fig:test5_RTA}
\end{figure}

The results obtained with the two collision models are very similar and consequently we only report those for the RTA approximation.
The computed results for the mass density at time $T_f=0.04$ with $V=5$ and $\gamma=0.002$ are given in figure \ref{fig:test5_RTA},
compared with a reference solution obtained with an explicit solver which is detailed in \cite{JP} with $N_x=400$. For the sake of clarity, we do not show the first order solution since it is very close to the second order one.

This problem ends with a steady solution which shows that our numerical method remains stable in stationary situations. Observe also that, as in the previous cases, the third order approximation gives the better results.



\section{Conclusions}

In this paper we have presented an Implicit-Explicit Runge-Kutta strategy for the numerical solution of the Boltzmann-Poisson system for
semiconductors in the diffusive scaling. We established sufficient conditions which permit to get high order in time asymptotic preserving and accurate schemes.

One of the main target of the present work was to overcome the severe time step restrictions to which this kind of problems are usually subject to.
A suitable recombination of the fluxes of the kinetic Boltzmann equation with those of the limit diffusive equation, together with an
IMEX strategy, allows us to achieve this goal and to make computation which are subject to a linear time step restriction $\Delta t = O(\Delta x)$.
These results have been obtained preserving high order accuracy in time and space for all regimes analyzed.
In order to consider more realistic models, a penalization technique which permits to avoid the costly inversion of non-linear collision
operators in the stiff regimes has also been considered. This allows to treat more realistic physical problems preserving efficiency and without loss of accuracy. In the last part we validated numerically the proposed approach. The results showed that the desired properties have been obtained.

In the future, we aim at treating more realistic simulations in the two (or three) dimensional setting and to perform a
stability analysis of the method developed.

\appendix
\section{Appendix}
\subsection{AP properties of the IMEX Runge-Kutta schemes}
\subsubsection{Analysis of type A schemes}
In this appendix we prove the AP property of the schemes proposed.
We start analyzing the IMEX schemes of type A. For these schemes multiplying equations (\ref{stRK}) and (\ref{stJK}) by $\e^2$ and
then setting $\e=0$ leads to
\begin{eqnarray}
 \Delta t A\,{\bf Q}(R)=0, \label{eq:equiR} \\
 \Delta t A\, {\bf g}(R,J)=0,  \label{eq:equiJ}
\end{eqnarray}
where ${\bf Q}(R) = ( Q(R^{(k)}) )_k $ and ${\bf g}(R,J) = ( g(R^{(k)},J^{(k)}) )_k$, for $k=1,\ldots,\nu$ are the
vectors containing the stages values.
Then, since for hypothesis $A$ is invertible, we obtain from (\ref{eq:equiR}) that ${\bf Q}(R)=0$, which implies
\be
R^{(k)} = P^{(k)}M, \hspace{5mm} k=1,\ldots,\nu.
\label{eq:Re0}
\ee
Moreover, from (\ref{eq:equiJ}) it follows that ${\bf g}(R,J)=0$, which implies
\be
J = -\frac{1}{\l} \Big( \bv\cdot\nabla_\bx R-E\cdot\nabla_\bv R \Big).
\label{eq:Je0}
\ee
In other words, at each stage the parities $r$ and $j$ are projected over the equilibrium state.\\
In the limit $\varepsilon\rightarrow 0$, replacing (\ref{eq:Re0}) and (\ref{eq:Je0}) in (\ref{Pkkinsvolta}) we get
\begin{equation}\label{PdiscSt}
  P = \r^n {\bf e} + \Delta t \wt A \nabla_\bx \Big( (1-\mu) D \nabla_\bx \cdot P + \eta\, E \cdot P \Big)
                  + \Delta t \wt A\, \wt G {\bf e} + \Delta t A \mu \nabla_\bx \Big( D \nabla_\bx \cdot P \Big).
\end{equation}
If $\mu(0)=1$ the above equation is reduced to an arbitrary IMEX Runge-Kutta scheme for the drift-diffusion equation in which the
diffusion term is discretized implicitly. We recall in fact that the matrix $A$ corresponds to a diagonally implicit Runge-Kutta method.
Let us notice, that at the discrete level the two space derivatives which are related to the transport parts
of equations (\ref{stRK}) and (\ref{stJK}) must be treated with the same numerical scheme otherwise perfect cancellation of the term
$(1-\mu)D\nabla_\bx (P)$ in equation (\ref{PdiscSt}) is not possible, even if $\mu(0)$=1.

We analyze now the numerical solution for the macroscopic density in the limit $\varepsilon\rightarrow 0$. This becomes
\begin{equation}\label{PdiscSol}
 \r^{n+1} = \r^n + \Delta t \, \wt w^{T} \nabla_\bx \Big( (1-\mu) D \nabla_\bx \cdot P + \eta \, E \cdot P \Big)
                 + \Delta t \wt G + \Delta t \, w^{T} \mu D \Delta_{\bx\bx} P.
\end{equation}
Again if the condition $\mu(0)=1$ is satisfied and the same numerical treatment for the space derivatives is used
perfect cancellation of $(1-\mu) D \nabla_\bx \cdot P$ in equation (\ref{PdiscSol}) is guaranteed.
This means that the diffusive term in the limit is discretized implicitly.

To conclude this part, let us observe that an additional requirement may be demanded to our method,
namely that in the limit $\varepsilon\rightarrow 0$ the distribution function is projected over the
equilibrium $r^{n+1}\rightarrow \rho^{n+1}M$.
To find conditions which guarantee this property to be satisfied, let us analyze the equation for
the parity $r^{n+1}$
\be
r^{n+1} = r^n + \Delta t \wt w^{T}{\bf f}_1(R, J) + \Delta t w^T \mu \bv \cdot \nabla_\bx \left( {\bv\over\l} \cdot \nabla_\bx R \right)
              + \Delta t w^T {1\over\e^2} Q(R).
\ee
Replacing $\Delta t \, Q(R)$ in the above equation with its expression obtained from formula (\ref{stRK}), i.e.
$\e^2 A^{-1}\left(R-r^{n}{\bf e}-\Delta t \widetilde{A} {\bf f}_1(R,J)-\Delta t A {\bf f}_2(R)\right)$, leads to
\be\label{eq:nsA}
r^{n+1} = r^n + \Delta t \widetilde w^{T} f_1(R, J) + \Delta t w^{T}f_2(R) + w^{T} A^{-1}
\left( R - r^{n} \overline{e} - \Delta t \widetilde{A} f_1(R,J) - \Delta t A f_2(R) \right).
\ee
So now, the conditions to satisfy in order to project the distribution function over the equilibrium
distribution are $1=w^{T} A^{-1} e$, $\wt w^{T} = w^{T} A^{-1} \wt A$ and $w^{T} A^{-1} R =\rho^{n+1}M$.
This last condition depends on the stage values vector $R$. Thus the only possibility to satisfy the requirement is that
the IMEX scheme is Globally Stiffly Accurate (GSA). In this case, in fact, we automatically have
\begin{equation}
 w^T A^{-1}=(0,\ldots,0,1)^T,\quad P^{\nu}M=\rho^{n+1} M
\end{equation}
We can conclude that if the IMEX scheme is of type A and GSA then $\lim_{\varepsilon\to 0} r^{n+1}=\rho^{n+1}M$.

\subsubsection{Analysis of type CK schemes}
The request that the matrix $A$ is invertible can be very restrictive when high order methods are demanded.
However, under additional hypothesis, we can obtain schemes which are asymptotic preserving, asymptotically accurate with
implicit treatment of the diffusive term even when the first row of the implicit Runge Kutta method contains only zeros.
In details, we can state that if the IMEX scheme is of type CK and GSA with for initial data which are ``consistent'' with the limit problem
then in the limit $\varepsilon\rightarrow 0$ the scheme (\ref{stRK})-(\ref{numsolj}) becomes an IMEX scheme for the drift
diffusion equation (\ref{televen}).
The request of consistent initial data means that we want the initial distribution function $f$ to be a perturbation of the
equilibrium distribution $M$ with the perturbation going to zero as $\varepsilon\rightarrow 0$.
In formulas, for the parities $r$ and $j$, the consistency of the initial data with the limit problem (\ref{televen}) reads
\be
r_0=\rho_0 M+g^{\varepsilon},
\qquad
j_0= - \bv \cdot \nabla_\bx r_0 + {q\over m} E \cdot \nabla_\bv r_0,
\qquad
\lim_{\varepsilon\to 0} g^{\varepsilon}=0.
\label{eq:cid}
\ee
Under these hypothesis, in the limit $\varepsilon\rightarrow 0$, in addition to the AP and AA property, the solution is projected over the equilibrium at each time step: $r^{n+1}=\rho^{n+1}M$.

To prove this, let us rewrite all vectors relative to the stages separating the first component from the remaining ones;
for ${\bf Q}$, for example, we have
${\bf Q}(R) = \left( Q(R^{(1)}), \wh Q(\wh R) \right)$, where the second component is a vector in $\RR^{\nu-1}$.
We also define ${\bf f}_1 = \left(f_1(R^{(1)},J^{(1)}),\wh f_1(\wh R,\wh J)\right)^{T}$, $w=(w_1,\wh{w})^T$ and similarly
for all other vectors. Now for this type of schemes, rewriting (\ref{stRK}), we obtain
\begin{eqnarray}
 R^{(1)} = r^{n}, &  & \nonumber \\
 \wh{R} = r^{n} \widehat e & + & \Delta t \, \wt a \, f_1(R^{(1)},J^{(1)}) + \Delta t \, \wh{\wt A} \, \wh f_1(\wh R, \wh J) \nonumber \\
                           & + & \Delta t \, a \, \left( \frac 1{\e^2} Q(R^{(1)}) + f_2(R^{(1)}) \right)
                               + \Delta t \, \wh A \, \left( \frac 1{\e^2} \wh Q(\wh R) + \wh f_2(\wh R) \right)
 \label{eq:RkCK}
\end{eqnarray}
while from (\ref{numsolr})
\begin{eqnarray}
 r^{n+1} = r^n & + & \Delta t \, \wt w_1 \, f_1(R^{(1)},J^{(1)}) + \Delta t\, \wh{\wt w}^T \,\wh f_1(\wh R, \wh J) \nonumber \\
               & + & \Delta t \, w_1 \, \left( \frac 1{\e^2} Q(R^{(1)}) + f_2(R^{(1)}) \right)
                 +   \Delta t \, \wh w^T \, \left( \frac 1{\e^2} \wh Q(\wh R) + \wh f_2(\wh R) \right).
 \label{eq:rnp1CK}
\end{eqnarray}
Multiplying (\ref{eq:RkCK})
by $\e^2$ and then imposing $\e=0$ we obtain
$$
\Delta t\, a\, Q(R^{(1)}) + \Delta t \, \wh A \, \wh Q(\wh R) = 0.
$$
Since $R^{(1)}=r^{n}=\r^n M$, it holds that $Q(R^{(1)})=0$ and so, being $\wh A$ invertible,
we have that $R^{(k)}=P^{(k)} M$ for $k=2,\ldots,\nu$, in other words at each stage the distribution function is projected over the equilibrium state.

Finally, multiplying (\ref{eq:RkCK}) by $\e^2$ leads to
\begin{equation}
 \Delta t \, \wh Q(\wh R) = \e^2 \wh A^{-1} \left[ \wh R - r^{n} \widehat e
 - \Delta t \, \wt a \, f_1(R^{(1)},J^{(1)}) - \Delta t \, \wh{\wt A} \, \wh f_1(\wh R, \wh J)
 - \Delta t \, a\, f_2(R^{(1)}) - \Delta t \, \wh A \, \wh f_2(\wh R) \right],
\end{equation}
and thus, substituting this equality in (\ref{eq:rnp1CK}) and using again the fact that $Q(R^{(1)})=0$, we have
\begin{eqnarray}
 r^{n+1} = r^n
 + \Delta t \, \wt w_1 \, f_1(R^{(1)},J^{(1)}) + \Delta t\, \wh{\wt w}^T \,\wh f_1(\wh R, \wh J)
 + \Delta t \, w_1 \, f_2(R^{(1)}) + \Delta t \, \wh w^T \, \wh f_2(\wh R) \nonumber \\
 + \wh w^T \, \wh A^{-1} \left( \wh R - r^{n} \widehat e
 - \Delta t \, \wt a \, f_1(R^{(1)},J^{(1)}) - \Delta t \, \wh{\wt A} \, \wh f_1(\wh R, \wh J)
 - \Delta t \, a\, f_2(R^{(1)}) - \Delta t \, \wh A \, \wh f_2(\wh R) \right). \label{eq:rnp1CK_l}
\end{eqnarray}
Again, to assure the projection over the equilibrium of the numerical solution $r^{n+1}=\rho^{n+1}M$, we need that $1=\wh w^T\wh A^{-1}\wh e$,
$\wt w_1=\wh w^T\wh A^{-1}\wt a$, $\wt{\wh w}^T = \wh w^T\wh A^{-1}\wt{\wh A}$,
$w_1=\wh w^T\wh A^{-1}a$, $\wh w^T=\wh w^T\wh A^{-1}\wh A$ and that $\wh w^T\wh A^{-1}\wh R = \rho^{n+1}M$. These last requirements are automatically satisfied if the IMEX scheme is also $GSA$.

\subsection{AP properties of the penalized IMEX Runge-Kutta schemes}
In this last part, we prove the AP property of the IMEX-RK schemes in the penalized case only for type A methods. A similar proof holds for the CK type scheme under the same additional hypothesis of the non penalized case.

Let us rewrite the IMEX scheme for (\ref{eq:penr}) in vector form.
The stages are given by
\begin{equation} \label{eq:penstR}
 R = r^{n} {\bf e} + \Delta t \widetilde{A} \, {\bf f}_1(R,J) + \frac{\Delta t }{\e^2} \wt A \Big( {\bf Q}(R)-{\bf L}(R) \Big)
   + \frac{\Delta t }{\e^2} A \, {\bf L}(R) + \Delta t A{\bf f}_2(R)
\end{equation}
while the numerical solution by
\begin{equation} \label{eq:pennsr}
 r^{n+1} = r^n + \Delta t \wt w^{T}{\bf f}_1(R, J) + \frac{\Delta t}{\e^2} \wt w^T \Big( {\bf Q}(R)-{\bf L}(R) \Big)
         + \frac{\Delta t}{\e^2} w^T {\bf L}(R) + \Delta t w^T {\bf f}_2(R),
\end{equation}
with ${\bf L}$ being the stages vector of the linearized collision term.
Multiplying (\ref{eq:penstR}) by $\e^2$ and imposing $\e=0$, we get the following equality
\begin{equation*}
 \Delta t \, \wt A \, \Big( {\bf Q}(R) - {\bf L}(R) \Big) + \Delta t  \, A \, {\bf L}(R) = 0,
 \end{equation*}
from which:
\begin{equation*}
 {\bf L}(R) = - A^{-1} \, \wt A \, \Big( {\bf Q}(R) - {\bf L}(R) \Big).
\end{equation*}
Since $A^{-1} \, \wt A$ is lower triangular with diagonal elements equal to zero, we get projection over the equilibrium at each stage
$$
{\bf L}(R^{(k)})=0\quad\Rightarrow\quad R^{(k)}=P^{(k)}M,\quad k=1,\ldots,\nu.
$$
Concerning the limiting numerical solution we isolate the term $\Delta t\,{\bf L}(R)$
in (\ref{eq:penstR}) and then we substitute it in (\ref{eq:pennsr}) obtaining
\begin{eqnarray}
  r^{n+1} &=& r^n + \Delta t \wt w^{T}{\bf f}_1(R, J)
            + \frac{\Delta t}{\e^2} \wt w^T \Big( {\bf Q}(R)-{\bf L}(R) \Big) + \Delta t w^T{\bf f}_2(R)  \nonumber  \\
          &+& w^T \, A^{-1} \left( R - r^{n} {\bf e} - \Delta t \, \wt A \, {\bf f}_1( R, J )
            - \frac{\Delta t}{\e^2} \wt A \Big( {\bf Q}(R)-{\bf L}(R) \Big) - \Delta t \, A \, {\bf f}_2(\wh R) \right). \;\; \label{eq:penlrs}
\end{eqnarray}
Observe that, in this case, the numerical solution still depends on $ 1 / \e^2$, unless we require that our IMEX scheme satisfies $ \wt w^T = w^T \, A^{-1} \wt A$. This requirement is necessary not only to guarantee the correct projection over the equilibrium but also to be able to actually compute the solution in the limit.
It is easy to verify that the GSA condition is a sufficient condition which permits to guarantee the above requirement but also
$ 1 = w^{T} A^{-1} e $, $\wt w^{T} = w^{T} A^{-1} \wt A$, $w^T = w^T A^{-1} A$, $ w^{T} A^{-1} R = \rho^{n+1} M $. In other words that the schemes are AP and that the solution is projected over the equilibrium distribution at each time step.

\subsection{Examples of second and third order IMEX schemes}

We report here the Butcher tableaux of the second and third order schemes used in our simulations. Namely the second order ARS(2,2,2) scheme \cite{ARS}
$$
   \begin{array}{c|ccc}
     0       &   0      &  0         &  0       \\
     \gamma  &  \gamma  &  0         &  0       \\
     1       &  \delta  &  1-\delta  &  0       \\
     \hline
             &  \delta  &  1-\delta  &  0
   \end{array}
\qquad\qquad
   \begin{array}{c|ccc}
     0       &  0       &  0         &  0       \\
     \gamma  &  0       &  \gamma    &  0       \\
     1       &  0       &  1-\gamma  &  \gamma  \\
     \hline
             &  0       &  1-\gamma  &  \gamma
   \end{array}
$$
with $\gamma = 1 -\sqrt 2 / 2$ and $\delta = 1 - 1/(2\gamma)$ and the third order IMEX BPR-(3,5,3)  scheme \cite{BPR12}
$$
   \begin{array}{c|ccccc}
     0    &  0      &  0     &  0     &  0     &  0    \\
     1    &  1      &  0     &  0     &  0     &  0    \\
     2/3  &  4/9    &  2/9   &  0     &  0     &  0    \\
     1    &  1/4    &  0     &  3/4   &  0     &  0    \\
     1    &  1/4    &  0     &  3/4   &  0     &  0    \\
     \hline
          &  1/4    &  0     &  3/4   &  0     &  0
   \end{array}
\qquad\qquad
 \begin{array}{c|ccccc}
     0    &  0      &  0     &  0     &  0     &  0    \\
     1    &  1/2    &  1/2   &  0     &  0     &  0    \\
     2/3  &  5/18   &  -1/9  &  1/2   &  0     &  0    \\
     1    &  1/2    &  0     &  0     &  1/2   &  0    \\
     1    &  1/4    &  0     &  3/4   & -1/2   &  1/2  \\
     \hline
          &  1/4    &  0     &  3/4   & -1/2   &  1/2
   \end{array}
$$

\section*{Acknowledgments}
The author G. Dimarco was supported by the French ANR project BOOST.


\end{document}